\documentclass[12pt]{article}

\usepackage{amsmath, epsfig, cite}
\usepackage{amssymb}
\usepackage{amsfonts}
\usepackage{latexsym}
\usepackage{float}
\usepackage{color}
\usepackage{mathrsfs}
\usepackage{mathtools}

\newtheorem{thm}{Theorem}[section]

\newtheorem{cor}[thm]{Corollary}

\numberwithin{equation}{section}

\newcommand{\qed}{{\hfill$\square$}\medskip}

\begin{document}
	\begin{center}
		{\large\bf  A $q$-Supercongruence Motivated by Higher-Order Generalized Lehmer-Euler Numbers   } 
	\end{center}
	\vskip 2mm \centerline{Wei-Wei Qi}
	
	\begin{center}
		{\footnotesize MOE-LCSM, School of Mathematics and Statistics, Hunan Normal University, Hunan 410081, P.R. China\\[5pt]
			{\tt wwqi2022@foxmail.com} \\[10pt]
		}
	\end{center}
	
	\vskip 0.7cm \noindent{\bf Abstract.}Certain generalization of Euler numbers was defined in 1935 by Lehmer using cubic roots of unity, as a natural generalization of Bernoulli and Euler numbers. In this paper, we define a new polynomial related to the higher-order generalized Lehmer-Euler numbers and determine its a $q$-supercongruence.
	
	\vskip 2mm
	\noindent{\it MR Subject Classifications}: 11A07, 11B68, 33D15	
	
	\section{Introduction} 
	
	For positive integer $r$ and $\alpha$, the \textit {generalized Lehmer-Euler numbers}  of order $\alpha$, which is denoted by $W_{r,n}^{(\alpha)}$ is defined by
	\begin{align*}
		\begin{aligned}
			\sum_{k=0}^{\infty} W_{r,n}^{(\alpha)}\frac{t^n}{n!}=\left(\frac{r}{\sum_{j=0}^{r-1}e^{\zeta^jt}}\right)^\alpha=\left(\sum_{l=0}^{\infty}\frac{t^{rl}}{(rl)!}\right),
		\end{aligned}
	\end{align*}
	where $\zeta=\zeta_r$ is the primitive r-th root of unity. When $r=3$ and $\alpha=1$, $W_n=W_{1,n}^{(1)}$ is the original Lehmer-Euler numbers, which was introduced and investigated  by D. H. Lehmer \cite{w--0w} in $1935$. Recently, Komatsu and Liu \cite{w--1w} showed that
	\begin{align*}
		\begin{aligned}
			W_{r,n}^{(\alpha)}=\sum_{k=0}^{n}{\alpha+k-1 \choose k}{\alpha+n \choose n-k}(-\frac{1}{r})^k\sum_{\substack{i_1 + i_2 + \cdots + i_r = k \\ i_1, i_2, \dots, i_r \geq 0}}\frac{k!}{i_1!i_2!\dots i_r!}\left(\sum_{j=1}^{r}\zeta^{j-1}i_j\right).
		\end{aligned}
	\end{align*}	
	If $r=2$, we will get an expression of the Euler number of hight-order \cite{w--2w}:	
	\begin{align*}
		\begin{aligned}
			E_n^{(\alpha)}=\sum_{k=0}^{n}{\alpha+k-1 \choose k}{\alpha+n \choose n-k}(-\frac{1}{2})^k\sum_{j=0}^{k}{k \choose j}(k-2j)^n.
		\end{aligned}
	\end{align*}		

The $q$-binomial coefficients are defined as
\begin{align}
	{n\brack k}={n\brack k}_q=\begin{cases}
		\displaystyle\frac{\left(q;q\right)_n}{\left(q;q\right)_k\left(q;q\right)_{n-k}} &\text{if $0\leqslant k\leqslant n$},\\[10pt]0 &\text{otherwise.}
	\end{cases}\label{aa-1}
\end{align}
where  \textit{$q$-shifted factorial} (\textit{$q$-Pochhammer symbol}) is defined $(a,q)_n=(1-a)(1-aq)\dots(1-aq^{n-1})$. Throughout the paper, the \textit{$q$-integer} $[n]_q$ is defined by $[n]=[n]_q=\frac{1-q^n}{1-q}=(1+q+\dots+q^{n-1})$.The $n$-th cyclotomic polynomial is given by
\begin{align*}
	\Phi_n(q)=\prod_{\substack{1\le k \le n\\[3pt](n,k)=1}}(q-\zeta^k),
\end{align*}
where $\zeta$ denotes a primitive $n$th root of unity.

For polynomials $A(q)$, $B(q)$, $C(q)$ $\in$ $\mathbb{Z}[q]$, the $q$-congruence
	$A(q)/B(q)\equiv 0 \pmod{C(q)}$
is understood as $A(q)$ is divisible by $C(q)$ and $B(q)$ is coprime with $C(q)$. In general,
\begin{align*}
	A(q)\equiv B(q) \pmod{C(q)} \quad  \iff \quad A(q)-B(q)\equiv 0 \pmod{C(q)}.
\end{align*}

In $2007$, Pan  \cite[Theorem 1.1 and Lemma 2.4]{w--3w} proved that
\begin{align*}
	\begin{aligned}
		2&\sum_{j=1}^{(p-1)/2}\frac{1}{[2j]}+2Q_p(2,q)-Q_p(2,q)^2[p]
		&\equiv \left(Q_p(2,q)(1-q)+\frac{(p^2-1)(1-q)^2}{8}\right)[p] \pmod{[p]_q^2},
	\end{aligned}
\end{align*}
and
\begin{align*}
	\begin{aligned}
		2&\sum_{j=1}^{p-1}\frac{(-1)^j}{[j]}
		&\equiv 2\sum_{j=1}^{(p-1)/2}\frac{1}{[2j]}-\frac{(p-1)(1-q)}{2}-\frac{(p^2-1)(1-q)^2}{24}[p] \pmod{[p]_q^2},
	\end{aligned}
\end{align*}
where  $p$ is prime and  the $q$-Fermat quotient is defined by
\begin{align*}
	Q_n(m,q)=\frac{(q^m,q^m)_{n-1}/(q;q)_{n-1}-1}{[n]}.
\end{align*} 
The Ap\'{e}ry number $A_n$ and $A_n'$ are given by
\begin{align*}
	A_n=\sum_{k=0}^{n}{n+k \choose k}^2{n \choose k}^2 \quad and \quad A_n'=\sum_{k=0}^{n}{n+k \choose k}^2{n \choose k}.
\end{align*}
In $2012$, Guo and Zeng\cite[(1.9)]{w--2w-1} established the following congruence to show two supercongruences involving the Ap\'{e}ry number
\begin{align}
	\sum_{k=0}^{n-1}{n+k \choose k}^2{n-1 \choose k}^2 \equiv 0 \pmod{n}. \label{0-1}
\end{align}
They also conjectured a $q$-analogue of \eqref{0-1} as follows:
\begin{align}
	\sum_{k=0}^{n-1}q^{(n-k)^2}{n+k \brack k}^2{n-1 \brack k}^2 \equiv q^{(n-1)^2}[n] \pmod{[p]_{q^{n/p}}^2}.\label{0-2}
\end{align}
In $2020$, Gu and Guo\cite[Theorem 1.3]{w--2w-2} gave a positive answer to \eqref{0-2} by establishing the following more general result:
\begin{align*}
	\sum_{k=0}^{n-1}q^{(n-k)^2}{n+k \brack k}^2{n-1 \brack k}^2 \equiv q[n] \pmod{\Phi_n(q)^2}. 
\end{align*}
In $2021$, Guo, Schlosser and Zudilin \cite[Conjecture 1]{w--2w-3-0} conjectured that
\begin{align}
	\begin{aligned}
	\sum_{k=0}^{n-1}&q^{r(n-k)^2+(r-1)k}{n+k \brack k}^{2r}{n-1 \brack k}^{2r} \\
	&\equiv q^{(r-1)n+1}[n]-\frac{r(2r-1)(n-1)^2q(1-q)^2}{4}[n]^3 \pmod{[n]\Phi_n(q)^3}. \label{0-3}
\end{aligned}
\end{align}
where $r$ is positive integer.  Which were  confirmed by Liu and Jiang\cite[Theorem 1.2]{w--2w-3}.

In the past few years, congruences involving $q$-binomial coefficients are an interesting project, which attracted many authors' attention. For more studies on $q$-congruences, one may consult \cite{w--2w-2}, \cite{w--2w-3} --\cite{w--3w} and so on.

Inspired by the above,	referring to the higher-order  generalized Lehmer-Euler numbers, we set a polynomial $M_n^*(\alpha)$ by $M_n^*(\alpha) = \sum_{k=0}^{n}{\alpha+k-1 \choose k}{\alpha+n \choose n-k}$, where   $\alpha$ and $n$  are positive integer.
In this paper, we main show a $q$-supercongruence of $M_n^*(\alpha)$.
	
\begin{thm}Let $n$ be a positive  odd  integer and  $\alpha$ be a positive even integer  with $\alpha \leq n$. Then
\begin{align}
	\begin{aligned}
	\sum_{k=0}^{n-1}&q^{k+1 \choose 2}{\alpha+k-1 \brack k}_q{\alpha+n-1 \brack n-1-k}_q\\
	&\equiv 2[n](1-q)+\frac{2q^\alpha[n]-[\alpha]}{[\alpha]}
	-2[n]Q_n(2,q)-2[n]\sum_{k=1}^{\alpha}\frac{(-1)^k}{[k]}
	  \pmod{\Phi_n(q)^2}. \label{t1}
	\end{aligned}
\end{align}
Furthermore, we have
\begin{align}
	\begin{aligned}
		\sum_{j=0}^{n-1}&q^j\sum_{k=0}^{j}q^{k+1 \choose 2}{\alpha+k-1 \brack k}_q{\alpha+n-1 \brack n-1-k}_q\\
		&\equiv -[n]-\frac{[\alpha]-[n]}{q^\alpha}-\frac{2[\alpha][n]}{q^\alpha}Q_n(2,q)-\frac{2[\alpha][n]}{q^\alpha}\sum_{k=1}^{\alpha}\frac{(-q)^k}{[k]} \pmod{\Phi_n(q)^2}. \label{t2}
	\end{aligned}
\end{align}	
\end{thm}
	
Recall that the Euler numbers \{$E_n$\}  is given by
\begin{align*}
	\frac{2}{e^x+e^{-x}}=\sum_{n=0}^{\infty}E_n\frac{x^n}{n!}\quad \left(\lvert x \rvert\textless\frac{\pi}{2}\right), \quad \sum_{k=0}^{n}{2n \choose 2k}E_{2k}=0 \quad (n=1,2,\cdots),\quad	E_0=1.
\end{align*}
In view of \cite[Lemma 2.2]{w--3w-1}(or \cite{w--3w-0}): Let $m, n \in \mathbb{N}$, we have
\begin{align}
	\sum_{k=1}^{n}(-1)^kk^m=\frac{(-1)^n}{2}\left(E_m(n+1)+(-1)^nE_m(0)\right).\label{t0}
\end{align}
Therefore, by Fermat's little theorem and \eqref{t0}, we deduce that
\begin{align}
	\sum_{k=1}^{\alpha}\frac{(-1)^k}{k}\equiv \frac{(-1)^\alpha}{2}\left(E_{p-2}(\alpha+1)+(-1)^\alpha E_{p-2}(0)\right)\pmod{p}.\label{t0-a}
\end{align}

  Letting $n=p$,  $q\rightarrow 1$ in \eqref{t1} and \eqref{t2}, then combining \eqref{t0-a}, we immediately obtain the following congruences.
\begin{cor} For prime $p$ and  positive even integer $\alpha$. Then
	\begin{align*}
		\begin{aligned}
			\sum_{k=0}^{p-1}{\alpha+k-1 \choose k}{\alpha+p-1 \choose p-1-k}
			\equiv -1-2pq_p(2)-p\left(E_{p-2}(0)-E_{p-2}(\alpha)\right)  \pmod{p^2},
		\end{aligned}
	\end{align*}
and
\begin{align*}
	\begin{aligned}
		\sum_{j=0}^{p-1}\sum_{k=0}^{j}{\alpha+k-1 \choose k}{\alpha+p-1 \choose p-1-k}
		\equiv -\alpha-2\alpha p q_p(2)-\alpha p\left(E_{p-2}(\alpha+1)+E_{p-2}(0)\right) \pmod{p^2}, 
	\end{aligned}
\end{align*}
where $q_p(2)=\frac{2^p-1}{p}$ is Fermat quotient.\\
\end{cor}

\section{Proof of Theorem 1.1}
In order to achieve the proof Theorem $1.1$, we need the following two $q$-harmonic series congruences. For odd integer $n$,
\begin{align}	
	\sum_{k=1}^{n-1}\frac{(-1)^k}{1-q^k}\equiv 2\sum_{k=1}^{\frac{n-1}{2}}\frac{1}{1-q^{2k}}-\frac{n-1}{2} \pmod{\Phi_n(q)},  \label{a1}
\end{align}	
and
\begin{align}	
	\sum_{k=1}^{\frac{n-1}{2}}\frac{1}{1-q^{2k}}\equiv -\frac{Q_n(2,q)}{1-q} \pmod{\Phi_n(q)}.   \label{a2}
\end{align}
In view of  \cite[Lemma 2.4]{w--3w} and  \cite[Theorem 1.1]{w--3w}, \eqref{a1} and \eqref{a2} can be easily proved by the same method as them. In fact, \eqref{a1} and \eqref{a2} are the weaker versions of \cite[Lemma 2.4]{w--3w} and \cite[Theorem 1.1]{w--3w}, respectively.\\ 

Using  \eqref{a1} and \eqref{a2}, we obtain
\begin{align}
	\begin{aligned}	
	\sum_{k=1}^{n-1}\frac{(-1)^k}{1-q^{k+\alpha}}&=(-1)^\alpha\left(\sum_{k=1}^{n-1}\frac{(-1)^k}{1-q^{k}}+\sum_{k=n}^{n+\alpha-1}\frac{(-1)^k}{1-q^{k}}-\sum_{k=1}^{\alpha}\frac{(-1)^k}{1-q^{k}}\right)\\
	&\equiv -2\sum_{k=1}^{\alpha}\frac{(-1)^k}{1-q^k}-\frac{2Q_n(2,q)}{(1-q)}-\frac{1}{1-q^n}-\frac{n-1}{2}+\frac{1}{1-q^\alpha} \pmod{\Phi_n(q)}. \label{a3}
   \end{aligned}
\end{align}
For $k \neq n-\alpha$, we have
\begin{align}
	\begin{aligned}
		{\alpha+n-1 \brack \alpha+k}_q&=\frac{(1-q^{n+\alpha-1})\dots(1-q^n)\dots(1-q^{n-k})}{(1-q)\dots(1-q^{\alpha+k})}\\
		&\equiv (1-q^n)\frac{(-1)^kq^{-{k+1 \choose 2}}}{(1-q^{\alpha+k}){\alpha+k-1\brack k}} \pmod{\Phi_n(q)^2}. \label{a4}
	\end{aligned}
\end{align}
By \eqref{a4}, we get
\begin{align}
	\begin{aligned}	
\sum_{k=1}^{n-1}&q^{k+1 \choose 2}{\alpha+k-1 \brack k}_q{\alpha+n-1 \brack n-1-k}_q-q^{n-\alpha+1 \choose 2}{n-1 \brack n-\alpha}_q{n+\alpha-1 \brack \alpha-1}_q\\
		&\equiv (1-q^n)\sum_{k=1}^{n-1}\frac{(-1)^k}{1-q^{k+\alpha}}+1 \pmod{\Phi_n(q)^2}. \label{a5}
	\end{aligned}
\end{align}
Furthermore, by the properties of the $q$-binomial coefficients, we deduce 
\begin{align}
	\begin{aligned}
	{n+\alpha-1 \brack \alpha-1}&={n-1 \brack n-\alpha}\frac{\prod_{i=1}^{\alpha-1}(1-q^{n+i})}{\prod_{i=1}^{\alpha-1}(1-q^{n-i})}\\
	&\equiv {n-1 \brack n-\alpha}(-1)^{\alpha-1} q^{{\alpha \choose 2}}\left(1+(1-q^n)\sum_{i=1}^{\alpha-1}\frac{1+q^i}{1-q^i}\right)\pmod{\Phi_n(q)^2}, \label{a6}
	\end{aligned}
\end{align}
and
\begin{align}
	\begin{aligned}
	{n-1 \brack \alpha-1}&=\frac{\prod_{i=1}^{\alpha-1}(1-q^{n-i})}{(q;q)_{\alpha-1}}
	&\equiv (-1)^{\alpha-1}q^{-{\alpha \choose 2}}\left(1-\sum_{i=1}^{\alpha-1}\frac{1}{1-q^i}\right)\pmod{\Phi_n(q)^2}, \label{a7}
	\end{aligned}
\end{align}
where we used the fact that $1-q^n\equiv 0\pmod{\Phi_n(q)}$.\\
Whence, combining \eqref{a6} and \eqref{a7} gives
\begin{align}
	\begin{aligned}
			{n-1 \brack \alpha-1}{n+\alpha-1 \brack \alpha-1}\equiv q^{-{\alpha \choose 2}}(\alpha-1)(1-q^n)-q^{-{\alpha \choose 2}} \pmod{\Phi_n(q)^2}. \label{a8}
	\end{aligned}
\end{align}

Since for any integer $t$,
\begin{align}
	\begin{aligned}
		q^{tn}&=1-(1-q^{tn})\\
		&=1-(1-q^n)(1+q^n+q^{2n}+\dots+q^{(t-1)n})\\
		&\equiv 1-t(1-q^n)\pmod{\Phi_n(q)^2}.\label{a9-0}
	\end{aligned}
\end{align}
It follows that
\begin{align}
	\begin{aligned}
		q^{{n-\alpha+1 \choose 2}-{\alpha \choose 2}}\equiv 1+\frac{2\alpha-n-1}{2}(1-q^n)\pmod{\Phi_n(q)^2}. \label{a9}
	\end{aligned}
\end{align}	
Substituting \eqref{a3}, \eqref{a8} and \eqref{a9} into \eqref{a5}, then simplifying we get
\begin{align}
	\begin{aligned}
		\sum_{k=1}^{n-1}q^{k+1 \choose 2}{\alpha+k-1 \brack k}_q{\alpha+n-1 \brack n-1-k}_q
		&\equiv 2[n](1-q)+\frac{[n](2q^\alpha-1)-[\alpha]}{[\alpha]}\\
		&-2[n]Q_n(2,q)-2[n]\sum_{k=1}^{\alpha}\frac{(-1)^k}{[k]}
		\pmod{\Phi_n(q)^2}, \label{a10}
	\end{aligned}
\end{align}
Noting that
\begin{align}
	\begin{aligned}
		\sum_{k=0}^{n-1}&q^{k+1 \choose 2}{\alpha+k-1 \brack k}{\alpha+n-1 \brack n-1-k}
		&=\sum_{k=1}^{n-1}q^{k+1 \choose 2}{\alpha+k-1 \brack k}{\alpha+n-1 \brack n-1-k}+{\alpha+n-1 \brack n-1}, \label{a11}
	\end{aligned}
\end{align}
and
\begin{align}
	\begin{aligned}
		{\alpha+n-1 \brack n-1}=(1-q^n)\frac{(1-q^{n+\alpha-1})(1-q^{n+1})}{(1-q)\dots(1-q^\alpha)}\equiv \frac{[n]}{[\alpha]} \pmod{\Phi_n(q)^2}.\label{a12}
	\end{aligned}
\end{align}
Then, combining \eqref{a10}--\eqref{a12}, we will arrive at \eqref{t1}.\\

Let us turn to \eqref{t2}. Similarly, we have
\begin{align}
	\begin{aligned}
		\sum_{j=0}^{n-1}&q^j\sum_{k=0}^{j}q^{k+1 \choose 2}{\alpha+k-1 \brack k}{\alpha+n-1 \brack n-1-k}\\
		&=\sum_{k=0}^{n-1}q^{k+1 \choose 2}{\alpha+k-1 \brack k}{\alpha+n-1 \brack n-1-k}\left([n]-[k]\right) \\
		&\equiv [n]-\sum_{k=1}^{n-1}q^{k+1 \choose 2}{\alpha+k-1 \brack k}{\alpha+n-1 \brack n-1-k}[k] \pmod{\Phi_n(q)^2},\label{b1}
	\end{aligned}
\end{align}
where we used \eqref{t1} and the fact $[0]=0$.\\
In addition, by \eqref{a4} we get
\begin{align}
	\begin{aligned}
	\sum_{k=1}^{n-1}&q^{k+1 \choose 2}{\alpha+k-1 \brack k}{\alpha+n-1 \brack n-1-k}[k]-[n-\alpha]{n-1 \brack n-\alpha}{\alpha+n-1 \brack n}q^{{n-\alpha+1}\choose 2}\\
	&\equiv (1-q^n)\sum_{k=1}^{n-1}\frac{(-1)^k[k]}{1-q^{k+\alpha}}+[n-\alpha]\\
	&=[n]\left(\sum_{k=1}^{n-1}\frac{(-1)^k}{1-q^{k+\alpha}}-\sum_{k=1}^{n-1}\frac{(-q)^k}{1-q^{k+\alpha}}\right)+[n-\alpha]
 \pmod{\Phi_n(q)^2}.\label{b2}
	\end{aligned}
\end{align}
Since $n$ is odd and $\alpha$ is even, using \eqref{a1} and \eqref{a2} yield
\begin{align}
	\begin{aligned}
	\sum_{k=1}^{n-1}\frac{(-q)^k}{1-q^k}&=\sum_{k=1}^{n-1}\frac{(-1)^k-((-1)^k-(-q)^k)}{1-q^k}=\sum_{k=1}^{n-1}\frac{(-1)^k}{1-q^k}-\sum_{k=1}^{n-1}(-1)^k\\
	&\equiv -\frac{n-1}{2}-\frac{2}{1-q}Q_n(2,q)\pmod{\Phi_n(q)}.\label{b3}
	\end{aligned}
\end{align}
Therefore, utilized \eqref{b3} we obtain
\begin{align}
	\begin{aligned}
		\sum_{k=1}^{n-1}\frac{(-q)^k}{1-q^{k+\alpha}}&=\frac{(-1)^\alpha}{q^\alpha}\left(\sum_{k=1}^{n-1}\frac{(-q)^k}{1-q^k}+(-q)^n\sum_{k=0}^{\alpha-1}\frac{(-q)^k}{1-q^{k+n}}-\sum_{k=1}^{\alpha}\frac{(-q)^k}{1-q^k}\right)\\
		&\equiv \frac{1}{q^\alpha}\left(\sum_{k=1}^{n-1}\frac{(-q)^k}{1-q^k}-2\sum_{k=1}^{\alpha}\frac{(-q)^k}{1-q^k}-\frac{q^n}{1-q^n}+\frac{q^{n+\alpha}}{1-q^{\alpha}}\right)\\
		&\equiv \frac{1}{q^\alpha(1-q)}\left(\frac{(1-n)(1-q)}{2}-2Q_n(2,q)-2\sum_{k=1}^{\alpha}\frac{(-q)^k}{[k]}-\frac{q^n}{[n]}+\frac{q^\alpha}{[\alpha]}\right)
		\pmod{\Phi_n(q)}.\label{b4}
	\end{aligned}
\end{align}
By \eqref{a9-0}, we have 
\begin{align}
	\begin{aligned}
	q^{{n-\alpha+1 \choose 2}-{\alpha+1 \choose 2}}\equiv \frac{1}{q^\alpha}+ \frac{(2\alpha-n-1)(1-q^n)}{2q^\alpha} 	\pmod{\Phi_n(q)^2}.\label{b5}
	\end{aligned}
\end{align}
With the help of  \eqref{a8} and \eqref{b5} gives
\begin{align}
	\begin{aligned}
	\left[n-\alpha\right]&{\alpha+n-1 \brack n}{n-1 \brack n-\alpha}q^{{n-\alpha+1}\choose 2}\\
	&\equiv \left(1+ (2\alpha-n-1)(1-q^n)\right)\frac{[\alpha]-[\alpha][n](\alpha-1)(1-q)-[n]}{q^\alpha}	\pmod{\Phi_n(q)^2} .\label{b6}
	\end{aligned}
\end{align}
Substituting \eqref{a3}, \eqref{b4} and \eqref{b6}  into \eqref{b2} and simplifying, we obtain
\begin{align}
	\begin{aligned}
		\sum_{k=1}^{n-1}&q^{{k+1 \choose 2}}{\alpha+k-1 \brack k}{n+\alpha-1 \brack n-1-k}[k]\\
		&\equiv \frac{[\alpha]-[n]}{q^\alpha}+\frac{2[\alpha][n]}{q^\alpha}Q_n(2,q)+\frac{2[\alpha][n]}{q^\alpha}\sum_{k=1}^{\alpha}\frac{(-q)^k}{[k]}	\pmod{\Phi_n(q)^2}.\label{b7}
	\end{aligned}
\end{align}
Finally, combining \eqref{b1} and \eqref{b7}, we get 
\begin{align*}
	\begin{aligned}
		\sum_{j=0}^{n-1}&q^j\sum_{k=0}^{j}q^{k+1 \choose 2}{\alpha+k-1 \brack k}{\alpha+n-1 \brack n-1-k}\\
		&\equiv -[n]-\frac{[\alpha]-[n]}{q^\alpha}-\frac{2[\alpha][n]}{q^\alpha}Q_n(2,q)-\frac{2[\alpha][n]}{q^\alpha}\sum_{k=1}^{\alpha}\frac{(-q)^k}{[k]} \pmod{\Phi_n(q)^2},
	\end{aligned}
\end{align*}
as desired.

Now, we complete the proof of Theorem $1.1$. \qed

\end{document}